\def\ove#1{\overline{#1}}    %
\def\ovs#1#2{\overset{#1}\to{#2}}
     \def\({\left(}       \def\al{\alpha}           
     \def\){\right)}      \def\e{\varepsilon}    
     \def\[{\left[}       \def\la{\lambda}
     \def\]{\right]}      \def\ffi{\varphi}

                                      \def\ot{\otimes}
     \def\<{\langle}                 \def\wh{\widehat}
     \def\>{\rangle}                 \def\wt{\widetilde}
                 \def\sbs{\subset}
\def\tr{\operatorname{trace}\,}

\def\id{\operatorname{id}}

\def\dim{\operatorname{dim}}

                     \def\sbs{\subset}

\def\span{\operatorname{span}}
        
\documentclass[12pt,oneside]{amsproc}
\usepackage[T2A]{fontenc}
\usepackage[cp1251]{inputenc}
\usepackage{amssymb}
\usepackage{amsthm}

\title{Some Thoughts on Approximation Properties}
\author{Oleg Reinov}
\address{ St. Petersburg State University,
Saint Petersburg, RUSSIA.}
\email{orein51@mail.ru}

\thanks{${ }$ Key words:   Approximation of operator;  approximation property; tensor product; eigenvalue.}

\begin{document}

                          $$ {} $$
\maketitle


\begin{abstract}
We study some known approximation properties and introduce and investigate
several new approximation properties, closely connected with different
quasi-normed tensor products. These are the properties like the $AP_s$
or $AP_{(s,w)}$ for $s\in (0,1],$ which give us the possibility to identify
the spaces of  $s$-nuclear and $(s,w)$-nuclear operators with the corresponding
tensor products (e.g., related to Lorentz sequence spaces). Some applications
are given (in particular,  we present not difficult proofs of
the  trace-formulas of Grothendieck-Lidskii type
for several ideals of nuclear operators).

 \end{abstract}

\vskip 0.3cm

Our main reference is [10]. All the notions , notations and facts, we use without any
reference, can be found in [1, 2, 4, 8, 10, 12].
\medskip

{\bf I.}\
The Grothendieck approximation property for a Banach space $X$ can be defined as follows:
$X$ has the $AP$ iff for every sequence $(x_n)_{n=1}^\infty\sbs X$ tending to zero, for any $\e>0$
there exists a finite rank (continuous) operator $R$ in $X$ such that  for each $n\in\Bbb N$
one has $||Rx_n-x_n||\le\e.$  Consider a natural question: for which sequences $(x_n)\in c_0(X),$
under some additional assumptions, the identity map $\id_X$ surely can be approximated
by finite rank operators, as above, and which of those conditions are sharp (or, if one wishes, optimal)?

One of the simplest fact (we think, known for more than 30 years) that

$(*)$\, if $(x_n)\in l_2(X),$ \,  $X$ is any, then the answer is positive.

Here ia a reason of this: Assuming $||x_n||\searrow 0,$ take any $N\in \Bbb N$ and consider
the linear span $\span[x_n]_1^N =: E_N$ as a subspace of $X.$ Define, fixing an $\e>0,$
a finite rank $R$ to be a projection from $X$ onto $E_N$ whose norm $\le \sqrt N.$

Now if $N$ is such that, for every $n\ge N,$ we have $||x_n||\le \frac {\e}{\sqrt N+1},$ then
$$
  ||Rx_n-x_n||=0 \ \text{  if  } \ n\le N,
$$
and
$$
 ||Rx_n-x_n||\le \e.
$$
Of course, instead of $(*)$ we can consider the the statement

$(**)$\, if $(x_n)\in l_{2,\infty}^0(X) $ \,  [Lorentz space with "o" small], $X$ is any, then
the answer is positive.

 The idea of the above proof is very simple and can be applied in some more general situations.
 For instance, every subspace of finite dimension of an $L_p$-space is
 $n^{|1/2-1/p|}$-complemented in that $L_p$-space. So,
  if $p\in [1,\infty],$  $\al=|\frac{1}{2}-\frac1p|$ and $X$ is a subspace of an $L_p$-space, then

 $(***)$\,  for every sequence $(x_n)\in l_{q,\infty}^0(X),$ where $1/q=\al,$
 the answer is positive.
 \smallskip

 {\it Remark 1}:\,
   About sharpness: it will be discussed a little bit later.
   \smallskip

   {\it Remark 2}:\,
   The statement $(***)$ has, as a matter of fact, the following  quantitative aspect:
   Given $\al\in [0,1/2]$ and a Banach space $X$ with the property that every
   finite dimensional subspace $F$ of $X$ is contained in a finite dimensional subspace $E\sbs X,$
   which $(E)$ in turn  is   $C\, (\dim F)^\al$-complemented in $X,$
   we have

   $(***)'$\,  for every sequence $(x_n)\in l_{q,\infty}^0(X),$ where $1/q=\al,$
   for any $\e>0$ there is a finite rank operator $R$ in $X$ so that
   $\sup_n ||Rx_n-x_n||\le\e.$

   Particular cases:

   \noindent
(i)\,
$q-2$ and $\al=1/2$ or $q=\infty$ and $\al=0$\,
($=$\, "$X$ is any Banach space" or "$X$ is isomorphic to a Hilbert space");

   \noindent
(ii)\,
$(x_n)\in l_q(X), q\in [2,\infty),$ or $(x_n)\in c_0(X), q=\infty$ [Hilbert case].
 \smallskip

 For a while let us introduce the notions of the corresponding approximation properties
 for a Banach space $X$ (taking into account that the possibility of approximations on $c_0$-sequences
 by finite rank operators gives us the
 Grothendieck's approximation property $AP):$
 Let $0<q\le \infty$ and $1/s=1/q+1.$
 We say that $X$ has the $\wt {AP}_s$  [resp., the $\wt{AP}_{s,\infty}]$
 if for every $(x_n)\in l_q(X)$ [resp., $l_{q,\infty}^0]$
 (where $l_q(X)$ means $c_0(X)$ for $q=\infty)$
 and for every $\e>0$ there exists a finite rank operator $R\in X^*\ot X$ such that
 $\sup_n ||Rx_n-x_n||\le \e.$
 Trivially, e.g., $\wt{AP}_{s_2} \implies \wt{AP}_{s_1}$ if $s_1\le s_2.$
 Thus, $\wt{AP}_1 (= AP)$ implies any $\wt{AP}_s.$

 The statement $(*)$ (and $(**)$) says that every Banach space has the above property
 $\wt{AP}_{2/3}$ (and even the $\wt{AP}_{2/3,\infty}).$ The statement $(***)$
 gives the corresponding result for $L_p$-subspaces. Moreover, the assertion mentioned
 in Remark 2, shows that, for instance, any subspace of any quotient (= any quotient of any
 subspace) of a Banach space of type 2 (resp., of cotype 2) and of cotype $p,$ \, $p\in[2,\infty)$
 (resp., of type $p'),$ possesses the $\wt{AP}_s$
 (even the $\wt{AP}_{s,\infty})$  with $1/s=1+|1/2-1/p|.$

 \vskip 0.23cm

{\bf  II.}\
Let us recall that the notion of the $AP$ of Grothendieck can be reformulated in terms
of the projective tensor products "$\wh\ot$". Namely, a Banach space $X$ has the $AP$
iff for every Banach space $Y$ the canonical (natural) mapping $Y^*\wh\ot X \to L(Y,X)$
is one-to-one (or, what is the same, the natural mapping $X^*\wh\ot X\to L(X):=L(X,X)$
is injective).
In [3], A. Grothendieck has considered also some other tensor products (linear subspaces   
of  "$\wh\ot$"'s), which we will denote by "$\wh\ot_s$" for $0<s\le1$ (so that $\wh\ot =\wh\ot_1):$
For Banach spaces $X$ and $Y,$ let $Y^*\wh\ot_s X$ be a subspace of the projective tensor product
$Y^*\wh\ot X$ consisting of the tensors $z\in Y^*\wh\ot X,$ which admit representations of the form
 \begin{equation}
z= \sum_{n=1}^\infty \la_n y'_n\ot x_n,     \label {(1)}           
\end{equation}
where $(\la_n)\in l_s,$\, $(y'_n)$ and $(x_n)$ are bounded sequences from $Y^*$ and $X$
respectively. With a natural "quasi-norm" (see [10])         
the linear subspace $Y^*\wh\ot_s X$ of the space $Y^*\wh\ot X$ can be considered as
a "quasi-normed tensor product" (it is then a complete metric space [3]).

One of the nice (with a non trivial proof in [3]) theorem of Grothendieck is the fact
that the natural map from $Y^*\wh\ot_{2/3} X$ into $L(Y,X)$ is injective for any Banach spaces $X,Y.$
Let us compare this Grothendieck's
result with a simple assumption in Section I, where "$s=2/3$" was appeared. Must be clear that it is
not a chance coincidence, and really we have
\smallskip

{\bf Theorem 2.1.}\,
For $s\in (0,1]$ and for a Banach space $X$ the following are equivalent:

$1)$\,
$X$ has the $\wt{AP}_{s}$ in the sense of the definition in Section I;

$2)$\,
$X$ has the $AP_s$ in the sense of the definition in [13], i.e.           
for every Banach space $Y$ the natural mapping $Y^*\wh\ot_s X\to L(Y,X)$ is one-to-one.
\smallskip

Let us mention also that
\smallskip

$(AP_s)$\,
A Banach space $X$ has the $AP_s,$ \, $0<s\le1,$ iff the canonical map
$X^*\wh\ot_s X \to L(X)$ is one-to-one (or, what is the same, there exists no tensor element
$z\in X^*\wh\ot_s X$ with $\tr z=1$ and $\wt z=0,$ where $\wt z$ is the associated (with $z)$
operator from $X$ to $X$).                                 
\smallskip

The analogous theorems and facts are maybe valid for the $\wt{AP}_{s,\infty}$ and the $AP_{s,\infty}$
from  [13] (see a small discussion below).                  
\smallskip

{\it Proof}\ of  the assertion $(AP_s).$\,
Suppose $X$ has the $AP_s,$ but there exists a Banach space $Y$ such that the natural map
$Y^*\wh\ot_s X\to L(Y,X)$ is not one-to-one.Take an element $z\in Y^*\wh\ot_s X$    which is
not zero, but generates a zero operator $\wt z: Y\to X.$
Then we can find an operator $U\in L(X,Y^{**})$ so that $\tr U\circ z=1.$
If $z=\sum_{k=1}^\infty \la_k\, y'_k\ot x_k$ is a representation of $z$ in $Y^*\wh\ot_s X$\,
($(\la_k)\in l_s, (x_k)$ and $(y'_k)$ are bounded), then
$$
 1=\tr z= \sum_{k=1}^\infty \la_k \<Ux_k, y'_k\>=\sum_{k=1}^\infty \la_k\, \<x_k, U^*y'_k\>
$$
and $\sum_{k=1}^\infty \la_k\, U^*y'_k(x) x_k=0$ for every $x\in X.$
Put $x'_k:= U^*y'_k,$ $z_0:= \sum_{k=1}^\infty x'_k\ot x_k\in X^*\wh\ot_s X.$ We have
$$
\tr z_0=1,\ \wt z_0\neq 0
$$
(by assumption on $X).$ Consider a  1-dimensional operator $R= x'\ot x$ in $X$
with the property that $\tr R\circ z_0>0.$ Then
$$
 0<\tr R\circ z_0 = \sum_{k=1}^\infty \la_k\, \<x'_k, x\> \<x', x_k>= \sum_{k=1}^\infty \<U^*y'_k, x\> \<x', x_k\>
 $$
 $$= \<\sum_{k=1}^\infty \la_k\, \<Ux, y'_k\> x_k, x\>= \<x', \sum_{k=1}^\infty \la_k\, U^*y'_k(x) x_k\>=0.
$$
\smallskip

{\it Proof}\ of Theorem 2.1.
We will use the assertion $(AP_s).$

$1) \implies 2).$\,
Let $z\in X^*\wh\ot_s X$ and $\tr z=1.$
Write $z= \sum \la_k\, x'_k\ot x_k,$ where the sequences $x'_k)$ and $(x_k)$ are bounded and $(\la_k)\in l_s,$
$\la_k\ge0,$ $(\la_k)$ is non-increasing.
Then
$$
 z=\sum_{k=1}^\infty (\la_k^s\, x'_k)\ot (\la_k^{1-s})
$$
(recall that $1/s=1+1/q;$ so $1-1/s=1/q).$
The sequence $(\la_k^{1-s}x_k)$ is in $l_q(X).$ By 1), for every $\e>0$ there exists
a finite rank operator $R\in X^*\ot X$ such that $||R(\la_k^{1-s}x_k)-\la_k^{1-s}||\le\e$
for each $k\in \Bbb N.$ It follows that, for this operator $R,$
$$
  |\tr (z- R\circ z)| = |\sum_{k=1}^\infty \<\la_k^s x'_k, \la_k^{1-s} x_k- R(\la_k^{1-s} x_k)\>|
  \le \sum_{k=1}^\infty \la_k^s ||x'_k||\cdot \e\le const\cdot \e.
$$
Hence, for small $\e>0$ we have that, for an operator $R\in X^*\ot X,$
$$
  |\tr R\circ z| \ge 1/2
$$
and therefore $z$ generates a non-zero operator $\wt z.$
\smallskip

Before consider a proof of the implication $2)\implies 1)$ we will make some additional remarks.
We collect the remarks in
\smallskip

{\bf Lemma 2.1.}\
Let $s\in (0,1],$ $q\in (0,\infty],$ $1/s=1+1/q.$
For $a:=(a_k)\in l_1$ and $b:=(b_k)\in l_q$ we have
 \begin{equation}
(\sum_{k=1}^\infty |a_k b_k|^s)^{1/s}
   \le \sum_{k=1}^\infty |a_k| \cdot   (\sum_{k=1}^\infty |b_k|^q)^{1/q}.  \label {(2)}
\end{equation}
Moreover,
$$
  ||a||_{l_1} = \sup_{||b||_{l_q}=1} (\sum_{k=1}^\infty |a_k b_k|^s)^{1/s}
$$
(if $q=\infty$, the evident changes have to be made in $(2)$).
\smallskip

{\it Proof}\ of Lemma 2.1.
We may consider the case where $q\in (0,\infty).$ Putting $p:= 1/s$
(then $1/p'=1-s=s/q$ and $sp'=q)$, we obtain
$$
    \sum_{k=1}^\infty |a_k b_k|^s \le
    (\sum_{k=1}^\infty |a_k|^{sp})^{1/p} \cdot   (\sum_{k=1}^\infty |b_k|^{sp'})^{1/p'}=
    (\sum_{k=1}^\infty |a_k|)^s \cdot   (\sum_{k=1}^\infty |b_k|^q)^{s/q}.
$$

For the second part:
Let $a=(a_k)\in l_1.$ Take $b_k:= \frac{|a_k|^{1/q}}{||a||_{l_1}^{1/q}}.$ Then
$\sum_{k=1}^\infty |b_k|^{q}= \sum_{k=1}^\infty \frac{|a_k|}{||a||_{l_1}}=1$ and
$$
  (\sum_{k=1}^\infty |a_k b_k|^s)^{1/s}=
  = (\sum_{k=1}^\infty \frac{|a_k|^{s/q}}{||a||_{l_1}^{s/q}}\, |a_k|^s)^{1/s}=
   (\sum_{k=1}^\infty \frac{|a_k|^{s/q+s}}{||a||_{l_1}^{s/q}})^{1/s}
$$
$$
=  (\sum_{k=1}^\infty \frac{|a_k|^{s(1+1/q)}}{||a||_{l_1}^{s/q}})^{1/s}=
 (\sum_{k=1}^\infty \frac{|a_k|}{||a||_{l_1}^{s/q}})^{1/s}=
 \frac {(\sum_{k=1}^\infty |a_k|)^{1/s}}{||a||_{l_1}^{1/q}}=
 (/\sum_{k=1}^\infty |a_k|)^{1/s-1/q}= ||a||_{l_1}.
$$
\smallskip

{\it Proof}\ of Theorem 2.1 (continuation).

$2) \implies 1).$
Let $X$ has the $AP_s,$ but does not have the $\wt{AP_s},$ $1/s=1+1/q.$
Then there is a sequence $(x_n)\in l_q(X)$ (if $q=\infty,$ we consider a sequence from $c_0(X)=l_\infty^0(X)$)
such that there exists an $\e>0$ with the property that for any finite rank operator
$R\in X^*\ot X$ one has $\sup_n ||Rx_n-x_n||>\e.$
Consider the space $C_0(K; X)$ for $K:= \{x_n\}_{n=1}^\infty \cup \{0\}.$
Every operator $U$ in $X$ can be considered as a continuous function on $K$
with values in $X$ by setting $f_U(k):= U(k)$ for $k\in K.$
In particular, for the identity map $\id$ in $X$ and for any $R\in X^*\ot X$ we have
$$
  ||f_{\id}-f_R||_{C_0(K;X)}\ge\e.
$$
The subset $\mathcal R:=\ove{\{f_R:\ R\in X^*\ot X\}}^{C_0(K;X)}$ of $C_0(K; X)$
is a closed linear subspace in $C_0(K; X).$ So, there exists an $X^*$-valued measure
 $\mu=(x'_k)_{k=1}^\infty\in C_0^*(K; X)= l_1(\{x_n\}_{n=1}^\infty)\cup \{0\}; X)$ such that
 $\mu|_{\mathcal R}=0$ and $\mu(f_{id})=1.$
 In other words, we can find a sequence $(x'_k)$ with $\sum_{k=1}^\infty ||x'_k||<\infty$
such that $\sum_{k=1}^\infty \<x'_k, x_k\>=1$ and $\sum_{k=1}^\infty \<x'_k, Rx_k\>=0$
for any $R\in X^*\ot X.$

Define a tensor element $z\in X^*\wh\ot X$ by $z:= \sum_{k=1}^\infty x'_k\ot x_k.$
Since $(x_k)\in l_q(X)$ and $(x'_k)\in l_1(X^*),$ we get from Lemma 2.1 that
$$
  (\sum_{k=1}^\infty ||x'_k||^s\, ||x_k||^s)^{1/s} \le
    \sum_{k=1}^\infty ||x'_k||\cdot (\sum_{k=1}^\infty ||x_k||^q)^{1/q}.
$$
Therefore, $z\in X^*\wh\ot_s X,$ $\tr z=\sum_{k=1}^\infty \<x'_k, x_k\>=1$ and
$\tr R\circ z=0$ for every $R\in X^*\ot X.$ This means that $X$ does not have the $AP_s.$
\smallskip

After Theorem 2.1 is proved, we can make a conclusion:
$AP_s= \wt{AP}_s$ for any $s\in (0,1].$

 \vskip 0.23cm

{\bf  III.}\
Now we are going to discuss some questions around the properties $\wt{AP}_{s,\infty}$
and $AP_{s,\infty}.$ The $\wt{AP}_{s,\infty}$ was defined above. Recall the definition
of  the $AP_{s,\infty}$ from, e.g., [13]:
We say that a Banach space $X$ has the $AP_{s,\infty},$ $0<s<1,$ if for every Banach space $Y$
the natural mapping $Y^*\wh\ot_{s\infty} X\to L(Y,X)$ is one-to-one, where
$$
  Y^*\wh\ot_{s\infty} X= \{z\in Y^*\wh\ot X\!:\  z=\sum_{k=1}^\infty \la_k y'_k\ot x_k,\ (x_k)\,
  \text{ and}\, (y'_k) \text{ are  bounded, }  (\la_k)\!\in l_{s\infty}^0\}\!.
$$
Let us consider the connections between the $AP_{s,\infty}$ and the $\wt{AP}_{s,\infty}.$
For a partial discussion of this we need a lemma, which follows from Lemma 2.1
by interpolation in Lorentz spaces.
\smallskip

{\bf Lemma 3.1.}\
Let $s\in (0,1), q\in (0,\infty), 1/s=1+1/q, r\in (0,\infty].$
If $a=(a_k)\in l_1,$ $b=(b_k)\in l_{qr},$ then $ab:=(a_kb_k)_{k=1}^\infty \in l_{sr}.$
In particular, for $a\in l_1$ and $b\in l_{q\infty}$ the sequence $ab$ is
in $l_{s\infty}$ (thus, evidently, in $l_{s\infty}^0).$
\smallskip

{\it Proof}\ of Lemma 3.1 consist of  the applications of Lemma 2.1 and the general
interpolation theorem for the multiplication operator $\wt a,$ induced by a fixed
sequence $a=(a_k)\in l_1:$
$\wt a$ maps $(b_k)$ to $(a_kb_k).$

Namely, fix $s\in (0,1), q\in (0,\infty)$ with $1/s=1+1/q.$
Take $s_1, s_2\in (0,1)$ and $q_1, q_2\in (0,\infty)$ so that for some $\theta\in (0,1)$
we have
$$
  \frac1q = (1-\theta)\frac1{q_1} + \frac1{q_2}, \ 0<\frac1{s_2}<\frac1s < \frac1{s_1}<\infty, \
    0<\frac1{q_2}<\frac1q < \frac1{q_1}<\infty, $$
    and
  $$  \frac1{s_1}= 1+ \frac1{q_1}, \  \frac1{s_2}= 1+ \frac1{q_2}.
$$
By Lemma 2.1, $\wt a$ maps $l_{q_1q_1}$ into $l_{s_1s_1}$ and
 $\wt a$ maps $l_{q_2q_2}$ into $l_{s_2s_2}.$
 Applying, e.g., Theorem 5.3.1 from [1]          
 or other results from the pages 113-114 in [1],       
we get that $\wt a$ maps $l_{qr}$ into $l_{sr},$ $0<r\le \infty$
(note that $1/s=1+1/q= 1+(1-\theta)/q_1+\theta/q_2 = (1-\theta) +\theta +(1-\theta)/q_1 +\theta/q_2=
(1-\theta)(1+1/q_1)+\theta (1+1/q_2)=(1-\theta)/s_1+\theta/s_2$).
\smallskip

{\it Remark} 3.1:\
As a matter of fact, $l_1\cdot l_{q\infty}=l_{s1}$ in Lemma 3.1. We need now only
the above inclusion.
\smallskip

Now let $t\in (0,1],$ $p\in (0,\infty],$ $r\in (0,\infty]$ and
consider a tensor product $\wh\ot_{t; p,r},$ defined in the following way:
For a couple of Banach spaces $X, Y$ the tensor product
$Y^* \wh\ot_{t; p,r} X$ consists of those elements $z$ of the projective tensor product
$Y^*\wh\ot X$ which admit representations of the type
$$
  z=\sum_{k=1}^\infty a_kb_k\, y'_k\ot x_k; \ (y'_k) \text{ and } (x_k) \text{ are bounded, }
  (a_k)\in l_t,\, (b_k)\in l_{pr}
$$
(recall that everywhere here we consider $l_{p\infty}^0$ in the case $r=\infty).$
\smallskip

{\it Remark} 3.2:\
As was noted in Remark 3.1, $l_1\cdot l_{q\infty}=l_{s1} (\sbs l_{s\infty}^0\sbs l_{s\infty}),$
where $0<s<1, 1/s=1+1/q,$ We have also
$$
 l_{s1}=l_1\cdot l_{q\infty}^0\ \text{and}\ l_1\cdot l_{q\infty}=l_1\cdot l_{q\infty}^0
$$
(so, for example, in the definition of $\wh\ot_{1; q, \infty}$
one can assume that $(a_k)\in l_1$ and $(b_k)\in l_{q \infty}^0).$
Indeed, if we use the equality $l_1\cdot l_{q \infty}=l_{s1},$ take
$d\in l_{s1}$ (assuming $d=d^*=(d^*_k)$).            
Then $\sum_{k=1}^\infty k^{1/s}\, d^*_k/k <\infty,$ i.e. $\sum_{k=1}^\infty k^{1/q}\, d^*_k <\infty.$
Let $\e=(\e_k)$ be  a scalar sequence such that $\e_k\searrow 0$ and
$\sum_{k=1}^\infty \e_k^{-1} d^*_k k^{-1/q}<\infty.$ Put
$$
 \al_k:= \frac{d^*_k k^{1/q}}{\e_k,} \ \beta_k:= \frac{\e_k}{k^{1/q}}.
$$
Then $\al:= (\al_k)\in l_1$ and $\beta:= (\beta_k)\in l^0_{q \infty}.$
So, $d=\al \beta\in l_1\cdot l_{q \infty}^0.$
Another way (not to use "$l_{s 1}$"):
Let $0<q<\infty,$ $\al\in l_1,$ $\beta\in l_{q \infty}$ (assuming, without loss of generality,
that $\beta=\beta^*).$ Consider a sequence $\e:= (\e_k)$ such that $\e_k\searrow 0$ and
$\frac{}{}(\al_k/\e_k)\in l_1.$ Put $\wt \al:= \al/\e=(\al_k/\e_k)$ and  $\wt beta:=\e\beta=(\e_k\beta_k).$
Then $\wt \al\in l_1,$ $\wt \beta\in l^0_{q \infty}$ and
$\al\beta= \wt\al \wt\beta\in l_1\cdot l^0_{q \infty}.$
\smallskip

Let us say that $X$ has the $AP_{t; p,r},$ if  for every Banach space $Y$ and for $t, p, r$ as above
the canonical mapping $Y^*\wh\ot_{t; p, r} X\to L(Y,X)$ is one-to-one.

By Lemma 3.1, if $s\in (0,1)$ and $1/s=1+1/q,$ then $\wh\ot_{1; q, \infty}\sbs \wh\ot_{s, \infty}.$
Therefore, we get
\smallskip

{\bf Corollary 3.1.}\
If $s\in (0,1)$ and $1/s=1+1/q,$ then $AP_{s, \infty} \implies AP_{1; q, \infty}.$
\smallskip

Evidently, also $AP_{s,\infty}\implies AP_s$ (for $s\in (0,1)$).
\smallskip

{\bf Theorem 3.2.}\
 Let  $s\in (0,1), q\in (0,\infty)$ and $1/s=1+1/q.$ If $X$ has the $AP_{1; q,\infty},$
 then $X$ has the $\wt{AP}_{s,\infty}.$
In particular, $AP_{s,\infty} \implies \wt{AP}_{s,\infty}.$
\smallskip

{\it Proof}.\
It is enough to repeat word for word      
the proof of the implication $2) \implies 1)$ of Theorem 2.1 ("continuation"),
just changing "$l_{q}(X)$" by "$l^0_{q, \infty}$" (no necessity to apply Lemma 2.1 or Lemma 3.1).
\smallskip

{\it Remark} 3.3.\
In this moment (when I am writing the text) I do not know whether the implication
"$\wt{AP}_{s,\infty} \implies AP_{s,\infty}$" is true, for Banach spaces. Of course, no questions
about the cases where $0<s\le 2/3$ (but the reason is only that every Banach space has
the $\wt{AP}_{2/3,\infty}$ and the $AP_{2/3,\infty}).$
\smallskip

Let $0<r<1$ and  $0< w\le \infty.$ or  $r=1$  and  $0< w\le 1.$
For Banach spaces $X, Y$ denote by $Y^*\wh\ot_{(r,w)} X$ the subset of
$Y^*\wh\ot X$ consisting of tensors $z$ such that
$$
  z=\sum_{k=1}^\infty \lambda_k\, y'_k\ot x_k,\ \text{where} \ (y'_k) \text{ and } (x_k) \text{ are bounded and }
  (\la_k)\in l_{rw}.
$$

As was noted in Remark 3.1, if  $s\in (0,1), q\in (0,\infty), 1/s=1+1/q,$ then
$l_1\cdot l_{q \infty}= l_{s1}$ (in the sense of the product in Lemma 3.1).
In general case, where $0< q_1, q_2, t_1, t_2 \le \infty,$
one has
\begin{equation}
 l_{q_1 t_1}\cdot l_{q_2 t_2}\ \text{provided that:}\ \frac1{q_1}+\frac1{q_2}=\frac1s\ \text{and}\
   \frac1{t_1}+\frac1{t_2}=\frac1t.  \label{(3)}
 \end{equation}
We can introduce a new definition of approximation properties, which are connected
with Lorentz sequence spaces, namely:
Let  $0<r<1$ and  $0< w\le \infty.$ or  $r=1$  and  $0< w\le 1.$
A Banach space $X$ has the $AP_{(r,w)},$  if for every Banach space $Y$
the natural map $Y^*\wh\ot_{(r,w)} X\to L(Y,X)$ is one-to-one.

It follows (from Remark 3.1 or from (3)) that $AP_{1; q, \infty}= AP_{(s, 1)}$ (for $s\in(0,1)$ and  $1/s=1+1/q)$
and, more generally,  $AP_{t;p,r} = AP_{(s,u)}$ for $1/t+1/p=1/s$ and $1/t+1/r=1/u$\, $(t\in(0,1]).$

Therefore, we have (for $s\in(0,1)$)
$$
 AP_{s,\infty} \implies AP_{(s,1)} \implies \wt{AP}_{s,\infty}.
$$
Moreover, taking into account the equality $\wh\ot_{1;q,\infty}= \wh\ot_{(s,1)}$
and applying the arguments from the proof of the implication
"$\wt{AP}_{s}\implies AP_s$" of Theorem 2.1, we easily get
\smallskip

{\bf Theorem 3.3.}\
$AP_{(s,1)}=\wt{AP}_{s,\infty}.$
\smallskip

{\it Proof}.\
As was mentioned above, $AP_{(s,1)}\implies \wt{AP}_{s,\infty}.$
Let $X$ has the $\wt{AP}_{s,\infty},$ i.e.
for every sequence $(x_n)\in l_{q,\infty}^0$ (where $1/s=1+1/q)$ and every $\e>0$
there exists a finite rank operator $R\in X^*\ot X$ such that
$\sup_n ||Rx_n-x_n||<\e.$ Since $AP_{(s,1)}= AP_{1; q,\infty},$ it is enough
to show that if $Y$ is a Banach space, $z\in Y^*\wh\ot_{1; q,\infty} X$ and $z\neq0,$
then the corresponding operator $\wt z: Y\to X$ is not zero too.

Let $z=\sum_{k=1}^\infty a_kb_k\, y'_k\ot x_k$ be a representation of $z$ with
$(x_k), (y'_k)$ bounded, $(a_k)\in l_1,$ $(b_k)\in l^0_{q\infty}$ and $b_k\searrow0.$
Then $(\wt x_k:=b_kx_k)\in l^0_{q \infty}$ and, for an $\e>0$ small enough (to be he chosen),
we can find an operator $R\in X^*\ot X$ with the property that
$\sup_n ||R\wt x_n-\wt x_n||\le \e.$
Since $z\neq0,$ we can find an operator  $V\in L(Y^*, X^*)$ such that
$\sum_{k=1}^\infty a_k\, \<Vy'_k, \wt x_k\>=1.$
Now, when $V$ is chosen, we have
$$
  1= \sum_{k=1}^\infty a_k\, \<Vy'_k, \wt x_k-R\wt x_k\> +\sum_{k=1}^\infty a_k\, \<Vy'_k, R\wt x_k\>
  $$
  $$\le \e\, ||(a_k)||_{l_1}\,  ||V||\cdot const + |\sum_{k=1}^\infty a_kb_k\, \<R^*Vy'_k,  x_k\>|,
$$
and, if $\e$ is small enough, we get for the finite rank operator
$R^*V: Y^*\to X^*$ that
$$
 |\tr z^t\circ (R^*V)|= |\tr (R^*V)\circ z^t|= |\sum_{k=1}^\infty a_kb_k\, \<R^*Vy'_k, x_k\>|>0.
$$
The last sum is the nuclear trace of the tensor element $\sum_{k=1}^\infty a_kb_k\, R^*Vy'_k\ot x_k,$
which is a composition $R\circ z_0$ of the finite rank operator $R$ and the tensor element
$\sum_{k=1}^\infty a_kb_k\, Vy'_k\ot x_k,$ that
belongs to the tensor product $X^*\wh\ot_{1;q,\infty} X.$
It follows that both $z_0$ and $z$ generate the non-zero operators $\wt z_0$ and $\wt z.$
\smallskip

{\it Remark} 3.4.\
Because of the equality  $\wh\ot_{1;q,\infty}= \wh\ot_{(s,1)},$ it follows from the proof of Theorem 3.3
that $X$ has the $AP_{(s,1)}$ iff the canonical mapping $X^*\wh\ot_{(s,1)}X \to L(X)$ is one-to-one
(just like in the case of the classical Grothendieck approximation property).
\smallskip

{\it Remark} 3.5.\
Of course, it follows from Theorem 3.3 that every Banach space has the $AP_{(2/3,1)},$ but
it is trivial because of the implication
$$
 AP^0_{(2/3,\infty)} \equiv AP_{2/3,\infty} \implies AP_{(2/3,w)} \ \text{for any } w<\infty
$$
(and, again, since every $X$ has the $AP_{2/3,\infty}!).$
\smallskip

Our question in Remark 3.3 can be reformulated now as:

$({}^*)$\ Is it true that the $AP_{(s,1)}$ implies the $AP_{s,\infty}?$

 \vskip 0.23cm

{\bf  IV.}\
Let us consider an application of the previous considerations.
Now we know, in particular, that every Banach space has the $AP_{(2/3,1)}.$
On the other hand, the corresponding operator ideal $N_{(2/3,1)}$
(related to the Lorentz space $l_{2/3\, 1 })$ has the eigenvalue type $l_{1}$
(see, e.g., [4, p. 243]). Since the continuous trace is unique on
$\wh\ot_{(2/3,1)}$ and $\wh\ot_{(2/3, 1)}= N_{(2/3, 1)},$
it follows from White's results [17] that for each Banach space $X$ and for every operator
$T\in N_{(2/3,1)}(X,X)$ the (nuclear) trace of $T$ is well defined and equals the sum of all
eigenvalues of $T:$
$$
  \tr T= \sum_{k=1}^\infty \mu_k(T)\ \text{(eigenvalues)}\, \forall\, X,\ \forall\, T\in N_{(2/3,1)}(X)
$$
(on the right is the so-called "spectral sum" of $T).$ More precisely, the last statement follows
from Theorem 4.1 below.

Let us explain in more details how we  apply a White's result. For this we formulate and
 prove a theorem which is almost immediate consequence of the White's theorem.
 \smallskip

 {\bf Theorem 4.1.}\
Let $A$ be a quasi-Banach operator ideal, $X$ be a Banach space, for which the set of all finite rank
operators is dense in the space $A(X).$ Suppose that the natural functional $\tr$ is bounded
on the subspace of all finite rank operators of $A(X)$ (and, therefore, can be extended to a continuous
functional "$\operatorname{trace}_A$" on the whole space $A(X)$).
  If the quasi-Banach operator ideal $A$ is of eigenvalue type $l_1,$ then
 the spectral trace (= "spectral sum") is continuous  on the space $A(X)$ and for any
 operator $T\in A(X)$ we have
 $$
   \operatorname{trace}_A(T) = \sum_{n=1}^{\infty} \mu_n(T).
 $$
 where $(\mu_n(T))_{n=1}^\infty$ is  the sequence
 of all eigenvalues of $T,$ counting by multiplicities.
 \smallskip

  {\it Proof}\ of Theorem 4.1.    \
  Let $T\in A(X).$
 By the assumption, the sequence $\{\mu_n(T)\}_{n=1}^\infty$ of all eigenvalues of $T,$ counting by multiplicities,
is in $l_1.$

Since the quasi-normed ideal $A$ is  of spectral (= eigenvalue)
type $l_1,$  we can apply the main result from the paper [17] of M.C. White, which asserts:  

$({}^*{}_*{}^*)$\,  {\it If $J$ is a quasi-Banach operator ideal with eigenvalue type $l_1,$ then
the spectral sum is a trace on that ideal $J$}.

Recall   (see  [12], 6.5.1.1, or Definition 2.1 in [17]) that a {\it trace}\ on an operator ideal $J$
is a class of complex-valued functions, all of which one writes as $\tau,$ one for each component
$J(E,E),$ where $E$ is a Banach space, so that

(i)\  $\tau(e'\otimes e)= \langle e',e\rangle$ for all $e'\in E^*, e\in E;$

(ii)\ $\tau(AU)=\tau(UA)$ for all Banach spaces $F$ and operators $U\in J(E,F) and    A\in L(F,E); $

(iii)\ $\tau(S+U)=\tau(S) +\tau(U)$ for all $S,U\in J(E,E);$

(iv)\ $\tau(\lambda U)= \lambda \tau(U)$ for all $\lambda\in \Bbb C$ and $U\in J(E,E).$

Our operator $T$  belongs to the space $A(X,X)=A(X)$ and
$A$ is of eigenvalue type $l_1.$ Thus, the assertion $({}^*{}_*{}^*)$ implies that
the spectral sum $\lambda,$ defined by
$\lambda(U):= \sum_{n=1}^\infty \lambda_n(U)$   for $U\in A(E,E),$
is a trace on $A.$

By principle of uniform boundedness (see [11], 3.4.6 (page 152), or [9]),
there exists a constant $C>0$ with the property that
$$
 |\lambda(U)|\le ||\{\lambda_n(U)\}||_{l_1} \le C\, a(U)
$$
for all Banach spaces $E$ and operators $U\in A(E,E). $

Now, remembering
that all operators in $A(X)$ can be approximated by finite rank operators and
taking in account the conditions (iii)--(iv) for $\tau=\lambda$,
we obtain that the $A$-trace, i.e. $\operatorname{trace}_A T,$
of our operator $T$ coincides
with $\lambda (T)$ (recall that the continuous trace is uniquely defined in
such a situation, that is on the space $A(X);$ cf. [12], 6.5.1.2).
\smallskip

 Since $\wh\ot_{1;2,\infty}= \wh\ot_{(2/3,1)},$ (see Theorem 3.3), we can reformulate the result,
 which we formulated in the very beginning of this section, as
 \smallskip

 {\bf Corollary 4.1.}\
 For each Banach space $X$ and for every operator
 $T\in N_{1;2,\infty}(X,X)$ the (nuclear) trace of $T$ is well defined and equals the sum of all
 eigenvalues of $T:$
 $$
   \tr T= \sum_{k=1}^\infty \mu_k(T)\ \text{(eigenvalues)}\, \forall\, X,\
   \forall\, T\in N_{(1;2,\infty)}(X).
 $$
 \smallskip

 {\it Remark \rm4.1}:\,
 Recall that A. Grothendieck [3] has obtained the assertion of the last fact for the
 case of $2/3$-nuclear operators, i.e. for the ideal
 $N_{2/3}= N_{(2/3,\, 2/3)}$ (note that $l_{2/3}\sbs l_{2/3\, 1}).$
 \smallskip

  \vskip 0.23cm

 {\bf  V.}\
 The discussion on Section I shows that, for $p\in[1,\infty],$
 any subspace of any quotient (= any quotient of any
subspace) of  an $L_p$-space possesses the $\wt{AP}_s$
 (even the $\wt{AP}_{s,\infty})$  with $1/s=1+|1/2-1/p|.$
 We apply now these facts together with the White theorem for proving
 some more theorems concerning the distributions of eigenvalues of the nuclear operators.
 Below we will use Theorem 2.1 and, therefore, the fact that  any subspace of any quotient
 of  an $L_p$-space possesses the ${AP}_s$ (where $p,s$ as above).
 Thus, for such Banach spaces $X,$ we can identify the tensor product $X^*\wh\ot_s X$
 with its canonical image in the space $L(X)=L(X,X),$ that is with the space
 $N_s(X)$ of all $s$-nuclear operators in $X,$ equipped with the quasi-norm induced from
 $X^*\wh\ot_s X.$

 We are going to give below the relatively simple proofs of some recent results from the papers
 [15] and [16]. Let us begin.
 \smallskip

 {\bf Theorem 5.1.}\
 Let $X$ be a subspace of an $L_p$-space,
 $1\le p\le \infty.$ If $T\in N_s(X,X),$\, where
 $1/s=1+|1/2-1/p|,$   \,
 then

 1.\, the (nuclear) trace  of $T$ is well defined,

 2.\, $\sum_{n=1}^\infty |\mu_n(T)|<\infty,$ where
 $\{\mu_n(T)\}$ is the system of all eigenvalues of the operator $T$
 (written in according to their algebraic multiplicities)

 and
 $$
  \operatorname{trace}\, T= \sum_{n=1}^\infty \mu_n(T).
 $$

  \vskip 0.2cm

 {\it Proof}.\
 Let $X$ be a subspace  of an $L_p$-space $L_p(\mu)$
 and $T\in N_s(X,X)$ with an s-nuclear representation
 $$
  T=\sum_{k=1}^\infty \la_k x'_k\otimes x_k,
 $$
 where $||x'_k||, ||x_k||=1$ and $\la_k\ge 0,$  $\sum_{k=1}^\infty \la_k^s<\infty.$
 By Hahn-Banach, we can find  the functionals $\wt x'_k\in L^*_p(/mu)$ \, $(k=1,2,\dots)$ with
 the same norms as the corresponding functionals $x'_k$ and so that
 $\wt x'_k|_X=x'_k$ for every $k.$
 Denote by $\wt T$  the operator
 $$
  \wt T: L_p(\mu)\to X,\ \wt T:= \sum_{k=1}^\infty \la_k \wt x'_k\otimes x_k,
 $$
and by $j: X\to L_p(mu)$ the natural injection. Since the space $X$ has the property $AP_s,$
we have $N_s(L_p(\mu), X)= L^*_p(\mu)\wh\ot_s X$ and, therefore,
the nuclear traces of the operators $j\wt T$ and $\wt Tj$ are well defined.
 We have a diagram
 $$
  X\ovs{j} L_p(\mu) \ovs{\wt T} X \ovs{j} L_p(\mu),
 $$
 in which $\wt Tj = T\in N_s(X).$
Hence, the complete systems of eigenvalues of the operators  $T=\wt Tj $
and $j\wt T\in N_s(L_p(\mu))$ coincide .
Applying Theorem 2.b.13 from [5] (see also [15]), we obtain that
the sequence $(\mu_k(j\wt T))$ is summable. Therefore, we  have
$\mu_k(T)\in l_1$
and we can apply Theorem 4.1. But we apply the theorem firstly for the simplest case
(later on we will continue the proof of our theorem 5.1).
 \smallskip

 The first assertion of the next theorem is due to A. Grothendieck [3],          
 the second one was proved by H,. K\"onig in [6]. Surprisingly,                
 but we could not find anywhere the main statement of the theorem about coincidence
 of the nuclear and spectral traces, neither in the monographs, nor in the mathematical journals.
 So we have no reference for this statement and have to formulate and to prove the next theorem
 here. Let us remark that, in any case, this theorem was proved (as a partial case of the proved there
 our Theorem 5.1) in [15].                              
 \smallskip

 {\bf Theorem 5.1'.}\
Let $L$ be  an $L_p$-space,
$1\le p\le \infty.$ If $T\in N_s(L, L),$\, where
$1/s=1+|1/2-1/p|,$   \,
then

1.\, the (nuclear) trace  of $T$ is well defined,

2.\, $\sum_{n=1}^\infty |\mu_n(T)|<\infty,$ where
$\{\mu_n(T)\}$ is the system of all eigenvalues of the operator $T$
(written in according to their algebraic multiplicities)

and
$$
 \operatorname{trace}\, T= \sum_{n=1}^\infty \mu_n(T).
$$

 \vskip 0.2cm

{\it Proof}.\
As was said above, the assertions 1 and 2 are well known.
To prove the last equality, consider the Banach operator ideal $\mathcal L_p$
of all operators which can be factored through $L_p$-spaces. Then the product
$\mathcal L_p\circ N_s$ is a quasi-Banach operator ideal of spectral (=eigenvalue)
type $l_1$ (e.g., by the assertion 2, proved earlier by H. K\"onig [6]).       
Now it is enough to apply Theorem 4.1 to finish the proof.
\smallskip

{\it Proof \, of Theorem }5.1 \, (continuation).
As was said, the complete systems of eigenvalues of the operators  $T=\wt Tj $
and $j\wt T\in N_s(L_p(\mu))$ coincide. By Theorem 5.1',
$$
\tr j\wt T = \sum_{k=1}^\infty \la_k \, \<\wt x'_k, jx_k\> = \sum_{n=1}^\infty \mu_n(j\wt T),
$$
the last sum is equal to
$$
 \sum_{n=1}^\infty \mu_n(T)
$$
 and the sum in the middle is
 $$
   \sum_{k=1}^\infty \la_k \, \<\wt x'_k, jx_k\> = \sum_{k=1}^\infty \la_k\, \<x'_k, x_k\>=\tr T.
 $$
 The (nuclear) trace of the operator $T$ is well defined, because the space $X$ has the $AP_s.$
 Therefore,
 $$
   \tr T=  \sum_{n=1}^\infty \mu_n(T),
 $$
 and we are done.
  \smallskip

 If $Y$ is a quotient of an $L_p$-space, then, for a compact operator                                                         
 $U\in L(E,E),$ the adjoint    $U^*$ is also a compact operator and these two  operators
 have the same eigenvalues $\mu\neq0$ with the same multiplicities
 (see, e.g., [11], Theorem 3.2.26, or [2], Exercise VII.5.35).                 
 Also, any quotient of an $L_p$-space has the $AP_s$ (where $p,s$ are as above).
 So, it follows immediately from the just proved Theorem 5.1
  \smallskip

  {\bf Corollary 5.1.}\
   Let $Y$ be a quotient  of an $L_p$-space,
 $1\le p\le \infty.$ If $T\in N_s(Y,Y),$\, where
 $1/s=1+|1/2-1/p|,$   \,
 then

 1.\, the (nuclear) trace  of $T$ is well defined,

 2.\, $\sum_{n=1}^\infty |\mu_n(T)|<\infty,$ where
 $\{\mu_n(T)\}$ is the system of all eigenvalues of the operator $T$
 (written in according to their algebraic multiplicities)

 and
 $$
  \operatorname{trace}\, T= \sum_{n=1}^\infty \mu_n(T).
 $$

  \vskip 0.2cm

  We used above some facts from the section I. After Theorem 5.1 and its consequence are proved,
  we are ready to present a simple prove of the corresponding result on the subspaces of quotients
  of the $L_p$-spaces (recall that, again, all such spaces have the $AP_s$ with $s$ and $p$  satisfying
  the same conditions).
  \smallskip

   {\bf Theorem 5.2.}\
  Let $W$ be  a quotient of a subspace (= a subspace of  a quotient) of
  an $L_p$-space,
  $1\le p\le \infty.$ If $T\in N_s(W,W),$\, where
  $1/s=1+|1/2-1/p|,$   \,
  then

  1.\, the (nuclear) trace  of $T$ is well defined,

  2.\, $\sum_{n=1}^\infty |\mu_n(T)|<\infty,$ where
  $\{\mu_n(T)\}$ is the system of all eigenvalues of the operator $T$
  (written in according to their algebraic multiplicities)

  and
  $$
   \operatorname{trace}\, T= \sum_{n=1}^\infty \mu_n(T).
  $$

   \vskip 0.2cm

  {\it Proof}.\
  Let $L_p(\mu)$ be an $L_p$-space. Take Banach subspaces
  $X_0\sbs X\sbs L_p(\mu)$ and consider the quotient $X/X_0.$
  If $T\in N_s(X/X_0, X/X_0)$ (=$(X/X_0)^*\wh\ot_s X/X_0),$ then
  $T$ admits a factorization of the type
  $$
    X/X_0 \ovs{A} c_0 \ovs{D} l_1 \ovs{B}  X/X_0,
  $$
  where $A,B$ are continuous and $D$ is a diagonal operator with a diagonal from $l_s.$

  Denoting by $\ffi: X\to X/X_0$ the factor map from  $X$ onto $X/X_0$ and
  taking a lifting $\Phi: l_1\to X$ for $B$ with $B=\ffi\Phi,$
  we obtain that the maps $\ffi \Phi DA: X/X_0\to X/X_0$ and $\Phi DA\ffi: X\to X$
  have the same eigenvalues $\mu\neq0$ with the same multiplicities:
  $$
    X\ovs{\ffi} X/X_0 \ovs{A} c_0 \ovs{D} l_1 \ovs{\Phi} X \ovs{\ffi}  X/X_0,
  $$
  The spaces $X$ and $X/X_0$ have the $AP_s.$ Therefore, we have (cf. the proof of Theorem 5.1)
  $$
    \tr \ffi \Phi DA = \tr \Phi DA\ffi.
  $$
  Since $X$ is a subspace of  $L_p(\mu),$ we have, by Theorem 5.1,
  $$
   \tr \Phi DA\ffi= \sum_{n=1}^\infty \mu_n(\Phi DA\ffi).
  $$
  Therefore,
  $$
    \tr T = \tr BDA = \tr \ffi \Phi DA = \sum_{n=1}^\infty \mu_n(\Phi DA\ffi) $$
    $$ =
     \sum_{n=1}^\infty \mu_n(\ffi\Phi DA) = \sum_{n=1}^\infty \mu_n(BDA) =\tr T.
  $$
  \smallskip


   \vskip 0.23cm

  {\bf  VI.}\
    As is well known, in the classical case of  the   Grothendieck approximation property $AP$
  if $X^*$ has the $AP,$ then the space $X$ also has this property. We will show now that the same
  is true for all approximation properties which are under consideration in this paper.

  Denote by $\wh\ot_\alpha$ any of the tensor product
     $\wh\ot_s,$    $\wh\ot_{s,\infty},$    $\wh\ot_{t; p,r}$,      $\wh\ot_{(r,w)}$
      with the parameters (see above), for which all those tensor products are the  linear subspaces
    of the projective tensor product $\wh\ot.$ Also, let us say that a Banach space $X$
   has the $AP_\al,$ if it is possesses the corresponding approximation property (i.e., $AP_s,$
   $AP_{s,\infty}$ etc.).
\smallskip

We need the following  auxiliary result which may be of its own interest (compare with Remark 3.4).
\smallskip

{\bf Proposition 6.1}
A Banach space $X$ has the $AP_\al$ iff  the canonical map $X^*\wh\ot_\al X\to L(X)$ is
one-to-one.
\smallskip

{\it Proof}.\
Suppose that the canonical map $X^*\wh\ot_\al X\to L(X)$ is one-to-one,
but there exists a Banach space $Y$ such that the natural map
$Y^*\wh\ot_\al X\to L(X)$ is not injective. Let $z\in Y^*\wh\ot_\al X\to L(X)$
be such that $z\neq0$ and the associated operator $\wt z$ is a 0-operator.
Then we can find an operator $V$ from $L(Y^*,X^*)$ (the dual space to the
projective tensor product $Y^*\wh\ot X)$ so that $\tr V\circ z^t=1,$
where, as usual, $z^t$ is the transposed tensor element, $z*t\in X\wh\ot Y^*.$
Since $V\circ z^t\in X\wh\ot X^*$ and   $\tr V\circ z^t=1,$ the tensor element
$(V\circ z^t)^t$ (which, evidently, belongs to $X^*\wh\ot_\al X)$
is not zero. On the other hand, the operator induced by this element must be
a 0-operator. Contradiction.
\smallskip

  \noindent
  {\bf Proposition 6.2.}\
With the above understanding, if the dual space $X^*$ has the $AP_\al,$ then
$X$ has the $AP_\al$  too.
   \medskip

  \noindent
  {\it Proof}.\,
  We use Proposition 6.1.
  As is known [5], the projective tensor product
   $Y^*\widehat\otimes Y$  is a Banach subspace of the  tensor product
    $Y^*\widehat\otimes Y^{**}.$
    The tensor product  $Y^*\widehat\otimes_\al Y $ is  a linear subspace of  $Y^*\widehat\otimes Y,$
    as well as  $Y^*\widehat\otimes_\al Y^{**}$ is a linear subspace of  $Y^*\widehat\otimes Y^{**}.$
    Therefore, the natural map  $Y^*\widehat\otimes_\al Y \to Y^*\widehat\otimes_\al Y^{**} $
    is one-to-one. Now if $Y^*$ has the $AP_\al,$ then the canonical map
     $Y^{**}\widehat\otimes_\al Y^* \to L(Y^*,Y^*)$ is one-to-one.
     Since we can identify the tensor product $Y^{**}\widehat\otimes_\al Y^*$ with
     the tensor product $Y^{*}\widehat\otimes_\al Y^{**}$ (because of the "symmetries"
      in the definitions of the corresponding tensor products), it follows that
     the natural map $Y^*\widehat\otimes_\al Y \to L(Y,Y)$ is one-to-one.
     Thus, if $Y^*$ has the $AP_\al,$ then $Y$ has the $AP_\al$ too.
     \smallskip


       \smallskip

     \noindent
     {\it Remark 2}:\   The inverse statement is not true. For example,
    if $s\in(2/3,1],$ then there exists a Banach space, possessing the Grothendieck
    approximation property, whose dual does not have the $AP_s$ (it is well known for the case
    where $s=1).$
    Moreover,  if $s\in(2/3,1],$ then we can find a Banach space $W$ such that
    $W$ has a Schauder basis and $W^*$ does not have the $AP_s.$
    Indeed, let $E$ be a separable reflexive Banach space without the $AP_s$ (see [7] or [8]).
      Let $ Z$ be a separable space such that $ Z^{**}$ has a basis
    and there exists a linear homomorphism $ \varphi$ from $ Z^{**}$
    onto $ E^*$
    so that the subspace $ \varphi^*(E)$ is complemented
    in $ Z^{***}$ and, moreover,
    $Z^{***}\cong \varphi^*(E)\oplus Z^*$
    (see [7, Proof of Corollary 1]).
    Put $W:= Z^{**}.$ This (second dual) space  $W$ has a Schauder basis and its dual  $W^*$ does not have the $AP_s.$

    \smallskip


    \vskip 0.23cm

   {\bf  VII.}\
   Let us  consider some more notions of the approximation properties associated with some other
   tensor products.
   For Banach spaces $X$ and $Y$ and $r\in (0,1], p\in [1,2],$ define  a quasi-norm  $|| \cdot||_{N_{[r,p]}}$
   on the tensor product  $X^*\ot Y$ by

    \begin{equation}\nonumber
   \|u\|_{N_{[r,p]}}:=\inf\left\{\|(x_{i}')_{i=1}^{n}\|_{\ell_{r}(X^*)}\cdot\|(y_{i})_{i=1}^{n}
   \|_{\ell_{p'}^{w}(Y)}:\ u=\sum_{i=1}^{n}x_{i}'\otimes y_{i}\right\}
   \end{equation}
 Here we denote, as usual, by $l_r(X^*)$ and $l_q^w(Y)$ the spaces of $r$-absolutely
 summable and weakly $q$-summable sequences, respectively.

 Denote by   $X^{*}\widehat\otimes_{[r,p]}Y$  the completion of the space $\left(X^{*}\otimes Y,\ \|\cdot\|_{N_{[r,p]}}\right).$
 We have a natural continuous injection
  $$j_{[r,p]}: X^{*}\widehat\otimes_{[r,p]}Y\rightarrow X^{*}\widehat\otimes Y$$ with $||j_{[r,p]}||\le1.$

  Every element   $u\in X^{*}\widehat\otimes_{[r,p]}Y$ has a representation of the type $u=\sum_{i=1}^{\infty}x_{i}'\otimes y_{i},$
where $(x_{i}')_{i=1}^{\infty}\in\ell_{r}(X^{*})$ and $(y_{i})_{i=1}^{\infty}\in \ell_{p'}^{w}(Y)$.
     Consider the natural mappings
     $$
             X^{*}\widehat\otimes_{[r,p]} Y   \overset{\widetilde{j}_{[r,p]}}\to  X^{*}\widehat\otimes Y  {\widetilde{j }}\to  L(X,Y).
      $$
      The image of the tensor product  $X^{*}\widehat\otimes_{[r,p]} Y $ under the composition
      $\widetilde{j}_{[r,p]}:=\widetilde{j}\circ \wt j_{[r,p]}$
      is denoted by  $N_{[r,p]}(X,Y)$. This is a quasi-Banach space of the $(r,p)$-nuclear operators (the quasi-norm is induced from the
      tensor produce $X^{*}\widehat\otimes_{[r,p]} Y).$
      It is not difficult to see that every operator $T\in N_{[r,p]}(X,Y)$ admit a factorization of the kind
      $$
        X \ovs{A} c_0 \ovs{D_r} l_1 \ovs{i} l_p \ovs{B} Y,
      $$
      where $A,B$ are compact, $i$ is the injection, $D_r$ is a diagonal operator with a diagonal from $l_r.$
        Also, every operator which can be factored in such a way is in $T\in N_{[r,p]}(X,Y).$
      \smallskip

      By the analogous way, we define the tensor product $X^{*}\widehat\otimes^{[r,p]}Y$ and the quasi-normed operator ideals
      $N^{[r,p]}(X,Y).$ Namely,
      $X^{*}\widehat\otimes^{[r,p]}Y$ is a linear subspace of the projective tensor product $X^*\wh\ot Y,$
      consisting of tensor elements $z$ which admit a representation
      $$
        u=\sum_{i=1}^{\infty}x_{i}'\otimes y_{i},
        $$
        where $(x_{i}')_{i=1}^{\infty}\in\ell_{p'}^w(X^{*})$ and $(x_{i})_{i=1}^{\infty}\in \ell_{r}(Y).$
        Its canonical image in $L(X,Y)$ is the quasi-normed space $N^{[r,p]}(X,Y).$
        It is not difficult to see that every operator $T\in N^{[r,p]}(X,Y)$ admit a factorization of the kind
        $$
          X \ovs{A} l_{p'} \ovs{D_r} c_0 \ovs{i} l_1 \ovs{B} Y,
        $$
        where $A,B$ are compact, $i$ is the injection, $D_r$ is a diagonal operator with a diagonal from $l_r.$
        Also, every operator which can be factored in such a way is in $T\in N^{[r,p]}(X,Y).$

        It is clear that $T^*\in  N_{[r,p]}(Y^*,X^*)$ implies $T\in N^{[r,p]}(X,Y)$   and
$T^*\in  N^{[r,p]}(Y^*,X^*)$ implies $T\in N_{[r,p]}(X,Y).$

  Now we can define the notions of the corresponding approximation properties by the usual way.
  We say that he space $X$ has the $AP_{[r,p]}$ (respectively, the $AP^{[r,p]})$ if
  for every Banach space $Y$ the natural mapping $ Y^{*}\widehat\otimes_{[r,p]} X\to L(Y,X)$
  (respectively, $ Y^{*}\widehat\otimes^{[r,p]} X\to L(Y,X)$)     is one-to-one.
  It can be seen that
  a Banach space $X$ has the $AP_{[r,p]}$ (or $AP^{[r,p]})$ iff  the canonical map $X^*\wh\ot_{[r,p]} X\to L(X)$
  (or $X^*\wh\ot^{[r,p]} X\to L(X)$) is
  one-to-one (the proof is essentially the same as the proof of Theorem 6.1).
  Also, if $X^*$ has the $AP_{[r,p]}$ (or $AP^{[r,p]})$   then
$X$ has the $AP^{[r,p]}$ (or $AP_{[r,p]})$ (the proof is the as in Theorem 6.2).
 \smallskip

{\bf  Theorem 7.1.}\
Let   $1/r-1/p=1/2.$ Every Banach space has the properties $AP_{[r,p]}$
and $AP^{[r,p]}.$
\smallskip

{\it Proof}.\
Suppose that $X\notin AP_{[r,p]}$  where $1/r-1/p=1/2.$ Let $z\in X^*\wh\ot_{[r,p]} X$ be a n element such that
    $\tr z=1, \tilde z=0.$
   Since $z=\sum x'_k\ot x_k,$ where  $(x'_k)\in l_r(X^*) $ and $ (x_k) $ is weakly $ p'$-summable, the operator $ \tilde z $ can be factored as
     $$
    \tilde z:\ X\overset{A}\to l_\infty \overset{\Delta}\to l_1 \overset{j}\to l_p \overset{V}\to X,
    $$
   where all the operators are continuous, , $ j $ is an injection, $ \Delta $ is a diagonal operator with a diagonal from $ l_r. $
Since $\tilde z=0,$ we have $V|_{j\Delta A(X)}=0.$ Consider $S:= j\Delta AV: l_p\to l_p.$
Evidently, $S^2=0$ и $\tr S=\tr z=1.$
    Since $S\in N_r(l_p,l_p),$  its nuclear trace equals the sum  of all its eigenvalues (see  Theorem 5.1' above).
    This contradicts the fact that $S^2=0.$
   \vskip0.1cm


    We are ready to apply the above results to the investigation of eigenvalues problems
    for $N_{[r,p]}$-  and $N^{[r,p]}$-operators. The first theorem below was proved in [16]
    by using Fredholm Theory. The same proof can be applied for the second theorem.
\smallskip

{\bf  Theorem 7.2.}\
Let   $1/r-1/p=1/2.$
For every Banach space $X$ and every operator $T\in N_{[r,p]}(X)$, $\tr(T)$  is well defined and
if $(\mu_{i})_{i=1}^{\infty}$ is a system of all eigenvalues of $T,$ then $(\mu_{i})_{i=1}^{\infty}\in l_1$ and
     \begin{equation}\nonumber
     \tr(T)=\sum_{i=1}^{\infty}\mu_{i}.
     \end{equation}
  \smallskip

{\bf  Theorem 7.3.}\
Let   $1/r-1/p=1/2.$
For every Banach space $X$ and every operator $T\in N^{[r,p]}(X)$, $\tr(T)$  is well defined and
if $(\mu_{i})_{i=1}^{\infty}$ is a system of all eigenvalues of $T,$ then $(\mu_{i})_{i=1}^{\infty}\in l_1$ and
     \begin{equation}\nonumber
     \tr(T)=\sum_{i=1}^{\infty}\mu_{i}.
     \end{equation}
  \smallskip

Both theorems can be proved by the analogues methods and
the proofs are almost the same as the proof of Theorem 5.2 (by using Theorem 7.1).
So we omit it here.
    \smallskip



   \vskip 0.23cm

  {\bf  VIII.}\
The next examples are taken from [16], where one can find the corresponding
proofs. They show that all the above positive results
concerning approximation properties and trace-formulas are sharp.

  \vskip0.1cm


  \vskip0.1cm

  {\bf Example 8.1.}\
  Let $r\in(2/3,1], p\in(1,2], 1/r-1/p=1/2.$
  There exist Banach spaces  $E$ and $V,$ $z_0\in E^*\wh\ot V, S\in L(V,E)$
  so that for every $p_0\in [1,p)$

  1)\, $z_0\in E^*\wh\ot_{[r,1]} V;$

  2)\, $V$ has a basis;

  3)\, $V$ is the space of type $p_0$ and of cotype $2;$

  4)\, $S\circ z_0\in E^*\wh\ot_{[r,p_0]} E;$

  5)\, $\tr S\circ z_0=1;$

  6)\, the corresponding operator  $\widetilde{S\circ z_0}$  is a 0-operator and,
  therefore, has no nonzero eigenvalues.
    \vskip0.1cm


  {\bf Example 8.2.}\
  Let $r\in(2/3,1), p\in(1,2], 1/r-1/p=1/2.$
 There exist Banach spaces $E$ and $V,$ $z_0\in E^*\wh\ot V, S\in L(V,E)$
  so that for every $\epsilon>0$

  1)\, $z_0\in E^*\wh\ot_{[r+\epsilon,1]} V;$

  2)\, $V$ has a basis;

  3)\, $S\circ z_0\in E^*\wh\ot_{[r+\epsilon,p]} E;$

  4)\, $\tr S\circ z_0=1;$

  5)\, the corresponding operator  $\widetilde{S\circ z_0}$  is a 0-operator and
 therefore, has no nonzero eigenvalues.
  \vskip0.1cm

  {\bf Example 8.3.}\
  Let $r\in(2/3,1]$, $p\in(1,2], 1/r-1/p=1/2.$
 There exist two separable Banach spaces $X$ and $Z$ so that

  (i)\, $Z^{**}$ has a basis;

  (ii)\, $\exists \, V\in X^*\wh\ot Z^{**}: \ V=\sum_{k=1}^\infty x'_k\ot z''_k; $\ $(x'_k)$  weakly
  $p'_0$-summable for each $p_0\in [1,p);$ $(z''_k)\in l_r(Z^{**});$

  (iii)\, $V(X)\subset Z; $  the operator $V$ is not nuclear as a map from $X$ into $Z.$

 Moreover, there exists an operator $U:Z^{**}\to Z$ such that

  $(\alpha)$\, $\pi_ZU\in N^{[r,p_0]}(Z^{**},Z^{**})=Z^{***}\wh\ot^{[r,p_0]} Z^{**},\ \forall \, p_0\in[1,p);$

  $(\beta)$\, $U$ is not nuclear as a map from $Z^{**}$ into $Z;$

  $(\gamma)$\, $\tr \pi_ZU=1;$

  $(\delta)$\, $\pi_ZU: Z^{**}\to Z^{**}$ has no nonzero eigenvalues.
  \vskip1.1cm






\end{document}